\documentclass[graybox]{svmult}

\usepackage{type1cm}        
\usepackage{makeidx}         
\usepackage{graphicx}        
\usepackage{multicol}        
\usepackage[bottom]{footmisc}

\usepackage{newtxtext}       %
\usepackage{newtxmath}       

\makeindex             
                       
 \newcommand{\tZ}{\mathcal{Z}}            
  \newcommand{\tx}{\boldsymbol{x}}  
  \newcommand{\R}{\mathbb{R}}                   
    
  \newcommand{\Ee}{\mathsf{E}}
  \newcommand{\tV}{\boldsymbol{V}}

\begin{document}

\title*{Functional $K$ Sample Problem via Multivariate Optimal Measure Transport-Based Permutation Test}
\titlerunning{Functional $K$ sample OMT permutation test}  
\author{Šárka Hudecová, Daniel Hlubinka and Zdeněk Hlávka}

\institute{Šárka Hudecová \at Charles University, Faculty of Mathematics and Physics,
Department of Probability and Mathematical Statistics,
Sokolovská 83, Prague, Czech Republic, \email{hudecova@karlin.mff.cuni.cz}
\and Daniel Hlubinka \at Charles University, Faculty of Mathematics and Physics,
Department of Probability and Mathematical Statistics,
Sokolovská 83, Prague, Czech Republic, \email{hlubinka@karlin.mff.cuni.cz}
\and Zdeněk Hlávka \at Charles University, Faculty of Mathematics and Physics,
Department of Probability and Mathematical Statistics,
Sokolovská 83, Prague, Czech Republic, \email{hlavka@karlin.mff.cuni.cz}
}
%
%
\index{Hudecová, Š.} 
\index{Hlubinka, D.} 
\index{Hlávka, Z.} 

\maketitle

\abstract{The null hypothesis of equality of distributions of functional data coming from $K$ samples is considered. The proposed test statistic is multivariate and its components are based on pairwise Cram\'{e}r von Mises comparisons of empirical characteristic functionals. The significance of the test statistic is evaluated via the novel multivariate permutation test, where the final single $p$-value is computed using the discrete optimal measure transport. The methodology is illustrated by real data on cumulative intraday returns of Bitcoin.}

\section{Introduction}
\label{sec:Intro}

Functional Data Analysis (FDA) is widely applied in various practical applications where the observed data can be represented as (possibly discretized) functions. This includes situations
where the variable of interest is observed over a continuous time domain or over a high-frequency discrete time domain; setups often arising in economics, finance, biomedicine, engineering, and many other fields.  FDA has become a very popular approach to analyzing high-dimensional datasets in past few decades, \cite{horvath.kokoszka, ramsay2013functional}.  One of the basic statistical problems in FDA is the functional $K$ sample comparison, see \cite{cuesta, cuevas04, gorecki15, zhang2013analysis} and further references therein. 
Most studies focus on the null hypothesis of equality of means, resulting in a functional version of analysis of variance (ANOVA). However, the equality of covariance operators has also been investigated, \cite{guo2019}. 
 In this contribution, we are concerned with a more strict null hypothesis that the whole underlying distributions in the $K$ samples are equal. 
 
 In many FDA tests, the resulting test statistic is univariate and its (asymptotic) distribution is rather complex.  
It is, therefore, common to apply a permutation test, because it requires minimal distributional assumptions and leads to reasonable results even in small samples.  
 The significance  is then evaluated using so called permutation distribution obtained by random rearrangements of the data. The corresponding $p$-value is computed as a proportion of permutation test statistics that are more extreme than the value computed from the original data. 
 
 Unfortunately, the permutation principle cannot be easily generalized for multivariate test statistics, 
mainly due to the lack of ordering in dimension $d\geq 2$.
 A standard tool here is to compute $p$-values for each component of the
vector test statistic (so called partial $p$-values) and then to combine them, \cite[Chapter 6]{pesarin}. 
Recently, 
\cite{HHH} proposed a multivariate permutation test based on discrete optimal measure transport, \cite{peyre2019,villani2021}, which is referred to in the following as the OMT permutation test. 

Let $X_{j,1},\dots,X_{j,n_j}$ be independent and identically distributed functional observations coming from a distribution with a characteristic functional $\varphi_j$,  $j=1,\dots,K$, and let the $K$ samples be independent. Consider the null hypothesis  of equality of distributions 
\begin{equation}\label{eq:H0}
\mathcal{H}_0: \varphi_1=\dots=\varphi_K
\end{equation}
against a general alternative. 
For $K=2$, 
\cite{HHP} propose a Cram\'er von Mises type test statistic for $\mathcal{H}_0$ 
based on the empirical characteristic functionals. 
For $K>2$, the authors recommend to reformulate $\mathcal{H}_0$ using  $d=\genfrac(){0pt}{1}{K}{2}$ pairwise comparisons, i.e. 
\begin{equation}\label{eq:H0-part}
\mathcal{H}_{0}^{(j,l)}: \varphi_j=\varphi_l, \quad \text{ for } 1 \leq j < l \leq K.
\end{equation}
The test statistic is then computed for each pair $(j,l$), $j<l$, leading  to $d$-dimensional vector test statistic $\boldsymbol{T}=(T_1,\dots,T_d)^\top$. Subsequently, the  permutation test is applied to the univariate statistic $\max_{1\leq j\leq d} T_j$. It is clear that some information can be lost by the transformation of 
the multivariate test statistic $\boldsymbol{T}$ to its maximum. 
Hence, 
in this contribution, we follow the latter approach for $K>2$ samples, but  the significance of $\boldsymbol{T}$ is evaluated using the OMT permutation test.  
This approach not only uses the information from all the elements of $\boldsymbol{T}$, but also allows for an interpretation of the contributions of the individual pairwise comparisons to the rejection of the null hypothesis.  
That is, if the null hypothesis $\mathcal{H}_0$ is rejected, the OMT permutation test detects the pairs $(j,l)$ that contribute the most to the rejection. The benefits of this procedure are illustrated on a real dataset on cumulative intraday returns of Bitcoin prices.

%
%

The paper is organized as follows. The considered multivariate test statistic $\boldsymbol{T}$  is introduced in Section~\ref{sec:2sample}. Section~\ref{sec:OMT} describes the OMT permutation test and discusses related issues. The real data application is provided in Section~\ref{sec:data}.

\section{Multivariate Test Statistic for the $K$ Sample Problem}\label{sec:2sample}

A functional random variable $X$ is a random element taking values in the Hilbert space $\mathcal{X} = \mathcal{L}_2[0,1]$ of square integrable functions on $[0,1]$ such that $\mathsf{E} \int_{0}^1 X^2(t)\mathrm{d} t<\infty$. The space $\mathcal{X}$ is
equipped with an inner product
$\langle u,v\rangle = \int_{0}^1 u(t)v(t)\mathrm{d} t$, $u,v\in \mathcal{X}$. 
A distribution of $X$ is fully described by a characteristic functional 
$
\varphi(w) = \Ee \exp\{\mathrm{i}\langle w,X\rangle\}$ for  $w\in\mathcal{X},
$
where $\mathrm{i}=\sqrt{-1}$.

Consider first the two sample problem, i.e. $K=2$. Then \eqref{eq:H0} reduces to $\mathcal{H}_{0}^{(1,2)}: \varphi_1=\varphi_2$.
\cite{HHP}  proposed a test statistic for $\mathcal{H}_{0}^{(1,2)}$ based on a comparison of the empirical characteristic functionals of the two samples 
\begin{equation}\label{eq:T}
S_{1,2} = \int_{\mathcal{X}}|\widehat{\varphi}_1(w) - \widehat{\varphi}_2(w)|^2 \mathrm{d} Q,
\end{equation}
where $\widehat{\varphi}_j$ is the empirical characteristic functional in the $j$-th sample,  defined as
\begin{equation}\label{eq:ECF}
\widehat{\varphi}_j(w) = \frac{1}{n_j} \sum_{i=1}^{n_j} \exp\left\{ \mathrm{i} \langle w,X_{j,i} \rangle \right\}, \quad j=1,2,
\end{equation}
and $Q$ is a centered Gaussian measure on $\mathcal{X}$ with a covariance operator with a kernel function $v:[0,1]^2\to\R$. 
In practice, the functions $X_{j,i}$, $i=1,\dots, n_j$, $j=1,\dots,K$, are recorded on a finite set of points $t_1<\dots<t_J$ from $[0,1]$, and therefore, the inner products in \eqref{eq:ECF} need to be replaced by a suitable Riemann approximation of the integral. Furthermore, it suffices to specify $v$ on the measurement points, via a matrix $\tV = \Big(v(t_i,t_j)\Big)_{i,j=1}^J$. 
Further details on the numerical computation of $S_{1,2}$ and discussions on the choice of $\tV$ can be found in \cite{HHP}. 

\smallskip

Consider now the general  problem with $K>2$, the null hypothesis $\mathcal{H}_0$ in \eqref{eq:H0} and the pairwise partial hypothesis in \eqref{eq:H0-part}. For each pair $(j,l)$, $1\leq j<l\leq K$, let $S_{j,l}$ be the two sample test statistic defined in \eqref{eq:T} for samples $j$ and $l$. It follows from \eqref{eq:T}  that large values of $S_{j,l}$ indicate that the distribution in samples $j$ and $l$ may differ. Hence, 
define the test statistic 
\[
\boldsymbol{T}=(S_{1,2},S_{1,3},\dots,S_{K-1,K})^\top = (T_1,\dots,T_d)^\top,
\]
i.e. $\boldsymbol{T}$ is a vector with elements $S_{j,l}$ ordered according to the
lexicographic ordering of $(j,l)$. Then $\boldsymbol{T}$ takes values in $[0,\infty)^d$ and large values of its components indicate violation of the corresponding pairwise null hypothesis. 

\begin{remark}
It is clear that the dimension $d$ of $\boldsymbol{T}$ growths rapidly with the number of samples $K$. Remark that the dimension can be reduced, if required, by considering only some of the pairwise comparisons. 
For instance, $\mathcal{H}_0$ holds if and only if $\varphi_1=\varphi_j$ for all $j=2,\dots,K$. 
Therefore, one can consider only $\mathcal{H}_0^{(1,j)}$ for $j=2,\dots,K$ and get $\boldsymbol{T}$ of dimension $K-1$. 
\end{remark}

\section{Permutation Test Based on the Optimal Measure Transport}\label{sec:OMT}

 Let 
\[
\tZ=(X_{1,1},\dots,X_{1,n_1},\dots,X_{K,1},\dots,X_{K,n_K})
\]
be the ordered list of the pooled data. Denote as $\boldsymbol{T}_0$ the $d$-dimensional test statistic $\boldsymbol{T}$ computed as described in Section~\ref{sec:2sample} for the original data $\tZ$.  A permutation version of $\boldsymbol{T}$ is obtained as follows:  The elements of $\tZ$ are randomly permuted leading to a new ordered list $\tZ^*=(Z_{1,1}^*,\dots,Z_{1,n_1}^*,\dots, Z_{K,1}^*, \dots, Z_{K,n_K}^*)$. The test statistic $\boldsymbol{T}$ is then computed for the corresponding $K$ samples, with the $j$-th sample being $Z^*_{j,1},\dots,Z_{j,n_j}^*$. 
This procedure is repeated independently $B$ times, leading to permutation test statistics $\boldsymbol{T}_1,\dots,\boldsymbol{T}_B$.

The OMT permutation test from \cite{HHH} is based on the discrete $L_2$ optimal measure transport of the 
  set $\mathcal{T}=\{\boldsymbol{T}_0,\boldsymbol{T}_1,\dots,\boldsymbol{T}_B\}$ of $B+1$ points in $\R^d$ to a specified grid set $\mathcal{G}$ of $B+1$ points in the unit ball $\{\tx\in\R^d:\|\tx\|\leq 1\}$, whose choice is discussed later.  
  Namely, the optimal mapping
  $F^*$ 
is defined as a bijection from $\mathcal{T}$ to  $\mathcal{G}$ that minimizes the quadratic loss  
$
 \sum_{i=0}^{B} \|F(\boldsymbol{T}_i) - \boldsymbol{T}_i\|^2.
$ 
Note that computation of the discrete OMT $F^*$ is a standard optimalization task, see \cite[Chapter~3]{peyre2019}. In \texttt{R}, \cite{R}, it can be computed using package \texttt{clue}, \cite{clue}.

The  permutation $p$-value  is  computed as the relative frequency of permutation statistics that are  ''more extreme" than $\boldsymbol{T}_0$ from the center-outward perspective (i.e. that are more distant from the center), that is
\[
\widehat{p} = \frac{1}{B+1} \left(1+\sum_{b=1}^{B}
  \boldsymbol{1}\Bigl\{\| F^*(\boldsymbol{T}_{b})\| \geq \|F^*(\boldsymbol{T}_{0})\|
  \Bigr\}\right),
  \]
  where $\boldsymbol{1}\{\cdot\}$ is the indicator function. If the grid set $\mathcal{G}$ is specified as \eqref{eq:grid} below, one can use also a $p$-value $\widetilde{p}=1-\|F^*(\boldsymbol{T}_0)\|$.

\medskip

Apart from the single final $p$-value, the OMT permutation test provides also an interpretation of the partial contributions of the individual pairwise comparison.
  The quantity $(1-\widetilde{p})^2$ can be interpreted as an overall  non-conformity score that measures the deviation of the data from the null hypothesis, and $\mathcal{H}_0$ is rejected on level $\alpha\in(0,1)$ if and only if the non-conformity score  exceeds $(1-\alpha)^2$. Let $F^*(\boldsymbol{T}_0)=\big(F_1(\boldsymbol{T}_0),\dots,F_d(\boldsymbol{T}_0)\big)^\top$. Then
\[
(1-\widetilde{p})^2 = \|F^*(\boldsymbol{T}_0) \|^2 = \sum_{j=1}^d  F^*_j(\boldsymbol{T}_0)^2, 
\]
so the value $F_j^{*}(\boldsymbol{T}_0)^2$ can be interpreted as an absolute contribution of the $j$-th component of $\boldsymbol{T}$ (corresponding to $j$-th pairwise comparison) to the rejection of the composite null hypothesis $\mathcal{H}_0$. 
Denote as  $D_i = F_j^{*}(\boldsymbol{T}_0)/{\| F^*(\boldsymbol{T}_0)\|}$. Then $\sum_{i=1}^d D_i^2=1$, and the vector $(D_1^2,\dots,D_d^2)^\top$ can be  interpreted as a vector of relative contributions of the individual pairwise comparisons.  See \cite{HHH} for examples and more details. 

\medskip

The computation of the OMT permutation $p$-value requires a choice of the grid~$\mathcal{G}$. 
In the current context, each component of $\boldsymbol{T}$ takes values in $\R^+$ with large values indicating violation of the corresponding partial pairwise hypothesis. Therefore, in view of Section 3.2 in \cite{HHH}, a suitable grid set $\mathcal{G}$ is a subset of $\{\tx: \tx \in[0,1]^d, \|\tx\|\leq 1\}$. 
In view of  \cite{HHH}, it is beneficial to specify $B$ such that $B+1=n_R\cdot n_S$ for positive integers $n_R$ and $n_S$ and compute the grid points $\mathcal{G}=\{\boldsymbol{g}_{ij}\}_{i,j=1}^{n_R,n_S}$ as 
\begin{equation}\label{eq:grid}
\boldsymbol{g}_{ij}=\frac{i}{n_R+1} \boldsymbol{s}_j\quad i=1,\dots,n_R, \ j=1,\dots,n_S,
\end{equation}
where $\boldsymbol{s}_1,\dots,\boldsymbol{s}_{n_S}$ are distinct directional vectors from the set $\mathcal{S}^+=\{\boldsymbol{x}\in[0,1]^d: \|\boldsymbol{x}\|=1\}$, distributed as uniformly as possible over $\mathcal{S}^+$. These directions can be obtained from a sequence  $\{\tx_i\}_{i=1}^{n_S}$  of   low discrepancy points in $[0,1]^{d-1}$ as
 $\boldsymbol{s}_j = \tau(\tx_j)$, where $\tau$ is a mapping from $[0,1]^{d-1}$ to  $\mathcal{S}^+$ such that
 $\tau(\boldsymbol{X})$ has a uniform distribution on  $\mathcal{S}^+$ whenever $\boldsymbol{X}$ has a uniform distribution in $[0,1]^{d-1}$. 
 The mapping $\tau$ can be derived analogously as in \cite[Section 1.5.3]{fang}, see also Section~\ref{sec:data} for the application with $d=3$. Remark that the minimal $p$-value $\widehat{p}$ that can be obtained by the OMT permutation test with  grid from \eqref{eq:grid} is $1/n_R$. For instance, if $n_R=20$, then $1/n_R=0.05$ and results with this $p$-value have to be interpreted as significant.

\section{Real Data Application}\label{sec:data}

As an illustration of the proposed methodology, we analyze the intraday Bitcoin prices from 2022. 
The data come from a larger database available at \texttt{kaggle.com} repository\footnote{\texttt{https://www.kaggle.com/datasets/jkraak/bitcoin-price-dataset}}. We consider intraday data with frequency 20 min (so there are $J=72$ observations each day) and calculate the logarithmic cumulative intraday returns (CIR).    

%

\begin{figure}[htbp]
\centering
\includegraphics[width=\textwidth]{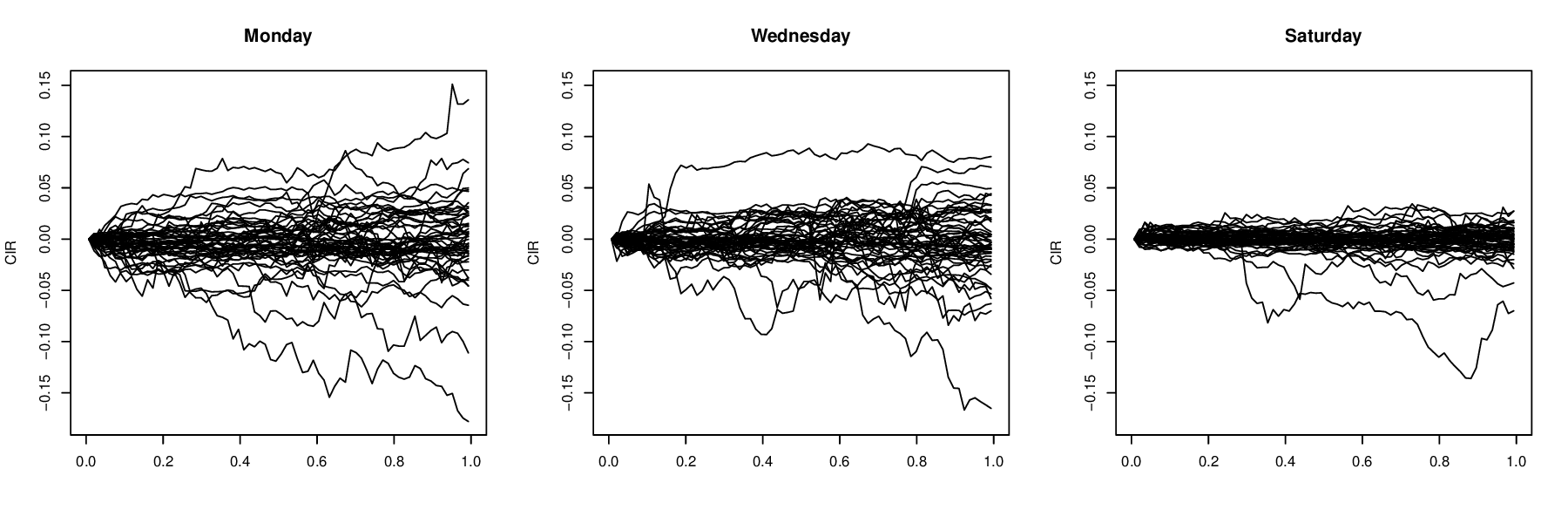}
\caption{Logarithmic cumulative intraday returns for Bitcoin in 2022: Mondays (52 curves), Wednesdays (52 curves), and Saturdays (53 curves).}\label{fig1}
\end{figure}

We are concerned with the question whether the behavior of CIR differs for working and weekend days. In particular, we consider the cumulative intraday returns for Mondays (sample 1), Wednesdays (sample 2) and Saturdays (sample 3), and test $\mathcal{H}_0$ for $K=3$.
The three samples, of sizes $n_1=52$, $n_2=52$, and $n_3=53$, respectively, are shown in Figure~\ref{fig1}.

Remark that the intraday Bitcoin prices for successive days are inherently dependent due to various market dynamics and external factors influencing price trends.  Moreover, they are unlikely to be identically distributed because of long-term trends and evolving patterns. The latter issue is   effectively addressed by transforming the data to CIR. The potential dependence between successive days is also partially mitigated by the CIR transformation. Furthermore, the selection of non-consecutive days (e.g., Monday, Wednesday, Saturday) reduces any lingering dependence between the observed data points. Consequently, any residual dependence between these chosen days can be reasonably assumed negligible.

The trivariate test statistic $\boldsymbol{T}$, whose components correspond to pairwise comparisons of samples (1,2), (1,3) and (2,3), respectively, was computed as described in Section~\ref{sec:2sample}. 
Based on empirical evidence, the authors from \cite{HHP} recommend to choose the matrix $\tV$ as data dependent, namely as an approximation to the inverse of the sample covariance matrix,  
because this choice leads to generally reasonable results under the null as well as against various alternatives. 
Following this recommendation, the following test statistics were considered:\\[1ex]
$-$ $\boldsymbol{T}^\text{inv}$ for $\tV$ being the approximate inverse to the sample covariance matrix computed from all data (without distinguishing the samples) and\\
$-$ $\boldsymbol{T}^\text{inv.pool}$ for $\tV$ being the approximate inverse to the pooled sample covariance matrix.

In both cases, the approximated inverse is computed from $9$ largest  eigenvalues. The permutation test statistics were calculated for $B=999$. 
The grid $\mathcal{G}$ in $\R^3$ is constructed as \eqref{eq:grid} with $n_R=40$ and $n_S=25$. 
The sequence $\{\tx_i\}_{i=1}^{n_S}$ is taken as a Halton sequence in $[0,1]^2$, see \cite{halton}, obtained in \texttt{R} using package \texttt{randtool}, \cite{randtool}. For $d=3$, the desired transformation from $\tx_j$ to $\boldsymbol{s}_j$ is via mapping $\tau=(\tau_1,\tau_2,\tau_3)^\top$ defined as
$
\tau_1(x_1,x_2) = 1-x_1$,
$\tau_2(x_1,x_2)= \sqrt{2x_1-x_1^2}\cos(\pi x_2/2)$, and $
\tau_3(x_1,x_2)= \sqrt{2x_1-x_1^2} \sin(\pi x_2/2)$.

\begin{figure}[htbp]
\centering
\includegraphics[width=0.48\textwidth]{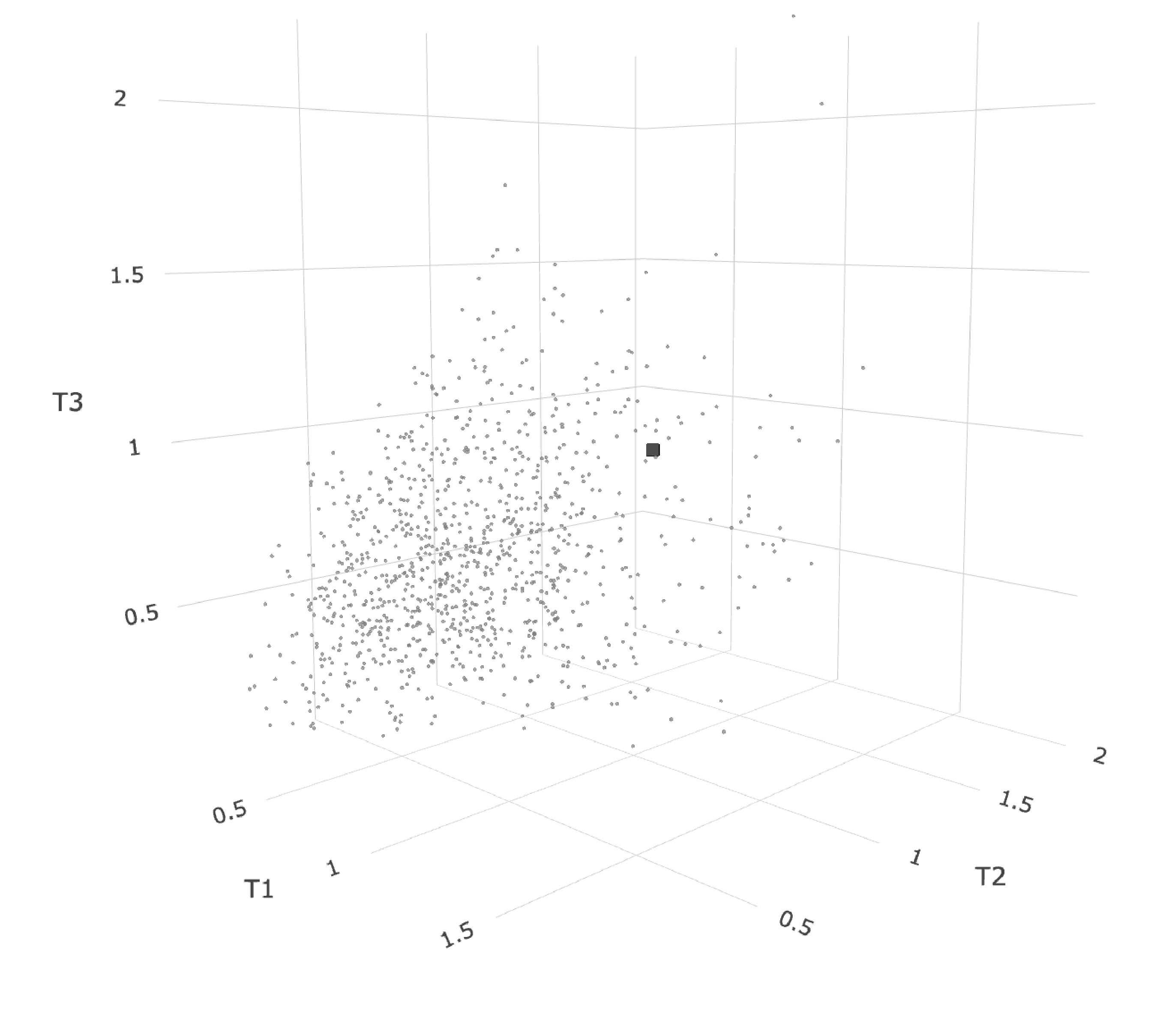}
\includegraphics[width=0.48\textwidth]{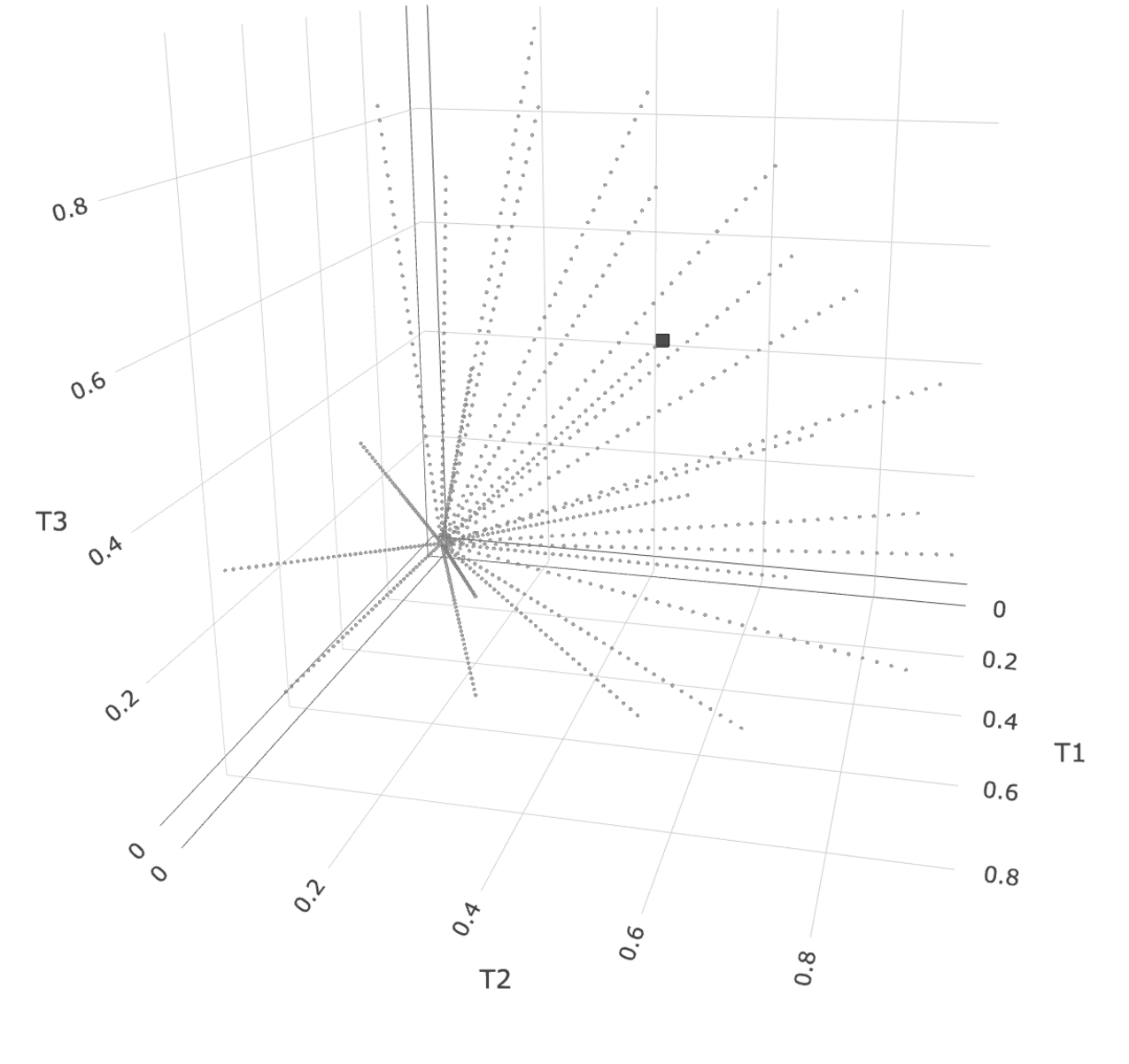}
\caption{Left panel: The test statistic $\boldsymbol{T}^\text{inv}$ (larger square) together with $B$ permutation replicas (small circles). All values multiplied by $10^4$. Right panel: The grid $\mathcal{G}$ with highlighted point $F^*(\boldsymbol{T}_0^{\text{inv}})$.}\label{fig2}
\end{figure}

 
 Figure~\ref{fig2} shows the test statistic  $\boldsymbol{T}^\text{inv}$ computed from the original data together with its $B$ permutation replicas, and the grid set $\mathcal{G}$ with the transformation $F^*(\boldsymbol{T}_0^{\text{inv}})$. It is visible that $F^*(\boldsymbol{T}_0^{\text{inv}})$, and subsequently $\boldsymbol{T}_0^{\text{inv}}$, is one of the most extreme points, so the obtained $p$-value is the minimal possible and significant. Namely, we get $\widehat{p}=0.025$ and  $\widetilde{p}=0.024$. The vector of individual contributions is $(D_1^2,D_2^2,D_3^2)=(0.353, 0.305, 0.343)^\top$, so it indicates that all three pairwise comparisons contribute equally to the rejection of the null hypothesis. 
 The conclusions for the test statistic $\boldsymbol{T}^{\text{inv.pool}}$ are very similar, so these are not discussed in detail.

 Table~\ref{tab1} presents $p$-values of various functional ANOVA tests for equality of mean functions, computed with the help of package \texttt{fdANOVA}, \cite{fdANOVA}. 
 In contrast to the results of our tests, all the tests for equality of mean function lead to highly non-significant $p$-values. This suggests that the violation of $\mathcal{H}_0$ is not due to the difference in means, but rather due to some other distributional aspects.  Figure~\ref{fig1} suggests that the difference in the three groups can be related to different covariance structures. Hence, we consider also the null hypothesis $\widetilde{\mathcal{H}}_0: C_1=C_2=C_3$, where $C_j$ is the covariance operator in sample $j$. This hypothesis can be reformulated in terms of three pairwise comparisons. For each pair, the operators can be compared using the test statistic based on the square-root distance, proposed by \cite{cabassi2017permutation}. This procedure results in  a trivariate test statistic that can be evaluated via the OMT permutation test. Using again $B=999$ and the same grid set, we get $\widehat{p}=0.05$ and $\widetilde{p}=0.049$. Note that due to the discreteness of possible $p$-values for the OMT permutation test, these values need to be interpreted as significant. The corresponding vector of individual contributions is $(0.062,0.909,0.028)^\top$, which reveals that $\widetilde{\mathcal{H}}_0$ is rejected mainly due to the difference in covariance operators of Mondays and Saturdays.

Finally note that our conclusions about differences in individual days are in agreement with some other empirical studies that examined differences in Bitcoin trading, see \cite{hansen2024} and references therein.

  \begin{table}[htbp]
  \centering
  \caption{Results for various functional ANOVA tests for equality of mean functions for the considered three samples. 
  FP is permutation test based on basis function representation from \cite{gorecki15}, CS is  $L_2$-norm-based parametric bootstrap test for  heteroscedastic samples from \cite{cuevas04}. L2B is $L_2$-norm-based test with  bias-reduced method of estimation, while L2b is $L_2$-norm-based bootstrap test, see \cite{zhang2013analysis}. FB is F-type test with   bias-reduced method of estimation, see \cite{shen,zhang2011}, and Fb is F-type bootstrap test, see \cite{zhang2013analysis}. 
  }\label{tab1}
\begin{tabular}{rrrrrrr}
  \hline\noalign{\smallskip}
Test & FP & CS & L2B & L2b & FB & Fb \\ 
$p$-value & 0.830 & 0.790 & 0.772 & 0.775 & 0.775 & 0.789 \\ 
  \noalign{\smallskip}\hline\noalign{\smallskip}
\end{tabular}
  \end{table}

\begin{acknowledgement}
The research was
supported by the Czech Science Foundation project GA\v{C}R No. 25-15844S.
\end{acknowledgement}

\bibliographystyle{abbrv}
\bibliography{HHH_lit.bib}

\end{document}